\documentclass[12pt]{article}
\usepackage{amsmath}
\usepackage{amssymb}
\usepackage{amsthm}
\usepackage{amscd}
\usepackage{amsfonts}
\usepackage{graphicx}%
\usepackage{fancyhdr}

\theoremstyle{plain} \numberwithin{equation}{section}
\newtheorem{theorem}{Theorem}[section]
\newtheorem{corollary}[theorem]{Corollary}

\newtheorem{lemma}[theorem]{Lemma}

\theoremstyle{definition}
\newtheorem{definition}[theorem]{Definition}

\newtheorem{remark}[theorem]{Remark}

\def\htm{H^2(\mu)}

\def\Scn{\mathcal{S}(\D^n)}

\def\={\ = \ }

\def\Om{\Omega}

\def\l{\lambda}
\def\z{\zeta}
\def\li{\lambda_i}
\def\lj{\lambda_j}          

\def\la{\langle}
\def\ra{\rangle}
\def\={ = }

\def\C{\mathbb C}

\def\T{\mathbb T}
\def\D{\mathbb D}

\def\be{\setcounter{equation}{\value{theorem}} \begin{equation}}
\def\ee{\end{equation} \addtocounter{theorem}{1}}
\def\beq{\begin{eqnarray*}}
\def\eeq{\end{eqnarray*}}

\def\bp{{\sc Proof: }}
\def\ep{{}{\hfill $\Box$} \vskip 5pt \par}

\def\bl{\begin{lemma}}
\def\el{\end{lemma}}
\def\bt{\begin{theorem}}
\def\et{\end{theorem}}
\def\bprop{\begin{prop}}
\def\eprop{\end{prop}}
\def\bd{\begin{definition}}
\def\ed{\end{definition}}
\def\br{\begin{remark}}
\def\er{\end{remark}}
\def\bexer{\begin{exercise}}
\def\eexer{\end{exercise}}

\newcommand{\g}{\gamma}
\newcommand{\lN}{\lambda_{N+1}}
\newcommand{\wN}{w_{N+1}}

\newcommand{\Frac}[2] {\displaystyle{\frac{#1}{#2}} }

\newcommand{\oo}{\omega}

\title{Hilbert function spaces and the Nevanlinna-Pick problem on the polydisc}
\author{David Scheinker}
\date{}                                           

\begin{document}

\maketitle

\begin{abstract}

In \cite{agmc_dv}, Agler and McCarthy used Hilbert function spaces to study the uniqueness properties of the Nevanlinna-Pick problem on the bidisc. 
In this work we give a geometric procedure for constructing a Nevanlinna-Pick problems on $\D^n$ with a specified set of uniqueness. On the way to establishing this procedure, we prove a result about Hilbert function spaces and partially answer a question posed by Agler and McCarthy.
\end{abstract}

\section{Introduction}

\subsection{Overview}
The Nevanlinna-Pick problem on $\D^n$ is to determine, given distinct nodes $\l_1,...,\l_N$ in $\D^n$ and target points $\oo_1,...,\oo_N$ in $\D$, whether there exists a function $F$ analytic on $\D^n$ with $||F||_\infty\leq1$ that satisfies $F(\li)=\oo_i$ for each $i$. A problem is called \textbf{extremal} if a solution $F$ satisfying $||F||_\infty=1$ exists and no solution $G$ satisfying $||G||_\infty<1$ exists.


Given a problem on $\D^n$, let $U$ denote the {\bf set of uniqueness} for the problem, the largest set on which all solutions agree. For $n=1$, if a problem is extremal, then Pick's 1916 results imply that the solution is unique, i.e. $U=\D$. On the other hand, for $n>1$, the following example shows that a problem may be extremal and 
yet fail to have a unique solution, i.e. $U\subsetneq\D^n$.\\

\noindent
\begin{exam}
For $n\geq2$, consider the problem with nodes\\
$(0,...,0), (1/2,...,1/2)\in\D^n$ and target points $0, 1/2\in\D$. 
Let $V=\{(z,...,z):z\in\C\}$. If $F$ is a solution, then $f(z)=F(z,...,z)$ is analytic on $\D$, satisfies $||f||_\infty\leq1$, $f(0)=0$, $f(1/2)=1/2$ and the classical 
Schwarz Lemma implies that $f(z)=z$. Thus, all solutions to the problem agree on $V\cap\D^n$, i.e. $V\cap\D^n\subset U$. Furthermore, that each coordinate function is a solution implies that $U=V\cap\D^n$.\\
\end{exam}

Various authors have studied the uniqueness properties of the Nevanlinna-Pick problem on the polydisc: in \cite{baltre98}, Ball and Trent show how to parameterize a certain class of solutions associated to a given Pick problem on $\D^2$; in \cite{agmc_three}, Agler and McCarthy classify those 2 and 3 point Pick problems on $\D^2$ that have a unique solution; in \cite{GHW} Guo, Huang and Wang give sufficient conditions for a 3 point Pick problem on $\D^3$ to have a unique solution; and in \cite{dsPLn}, a work closely related to this work, the present author gives a geometric procedure for constructing a Pick problem on $\D^n$ with a unique solution. 

The starting points of this work is a result from \cite{agmc_dv}. To state it, we say that an algebraic variety $V\subset \C^n$ is \textbf{inner} if each of it's irreducible components $V_i$ meets $\D^n$ and exits $\D^n$ through the n-torus, i.e. $V_i\cap\D^n\neq\emptyset$ and $V_i \cap \partial ( \D^n) \subset \T^n$. 

\bt{(Agler and McCarthy, \cite{agmc_dv})}
\label{agmc_ext}
Given an extremal Nevanlinna-Pick problem on $\D^2$, there exists a 1-dimensional inner variety $V$ such that all solutions to the problem agree on 
$V\cap\D^2$, i.e. $V\cap\D^2\subset U$.
\et


Theorem \ref{agmc_ext} suggests that one ask if every 1-dimensional inner variety $V\subset\D^2$ is contained in the uniqueness set of some problem. Our first result answers this question in the affirmative on $\D^n$.

\bt
\label{EachInner}
Given a 1-dimensional inner variety $V\subset\C^n$, there exists an extremal Nevanlinna-Pick problem on $\D^n$ such that $V\cap\D^n$ is the set of 
uniqueness for the problem, i.e. $V\cap\D^n=U$.
\et

Theorem \ref{EachInner} suggests that one ask if every inner variety $V\subset\C^n$, not necessarily 1-dimensional, is the set of uniqueness for some Pick 
problem on $\D^n$. Our second result gives a partial answer.

\bt
\label{kUniqueSet}
Given $k\leq n$ there exists a $k$-dimensional inner variety $V\subset\C^n$ and a Nevanlinna-Pick problem on $\D^n$ such that $V\cap\D^n$ is the 
set of uniqueness for the problem, i.e. $V\cap\D^n=U$.
\et
We prove Theorems \ref{EachInner} and \ref{kUniqueSet} in section 6.

In the remainder of this introduction, we discuss three results of independent interest that will be used to prove Theorems \ref{EachInner} and \ref{kUniqueSet}, Theorems \ref{FunctionSpace}, \ref{ThmB} and \ref{ThmD}.

\subsection{Theorem \ref{FunctionSpace}}
Theorem \ref{FunctionSpace} gives sufficient conditions for a Nevanlinna-Pick problem on a Hilbert function space to have a unique solution. Given a 1-dimensional inner variety $V$, we will construct a Hilbert function space on $V\cap\D^n$ and use Theorem \ref{FunctionSpace} to prove that for a certain type of problem on 
$\D^n$, all solutions agree on $V\cap\D^n$. To state the theorem, we need some notation and the notion of a ``Pick uniqueness kernel," defined below.

The multiplier algebra of a Hilbert function space $H(X)$, $Mult(H(X))$, is the normed algebra of functions $\phi$ on $X$ that satisfy $\phi f\in H(X)$ for each 
$f\in H(X)$. The norm is given by 
$||\phi||=||M_\phi||$, where $M_\phi$ is the bounded linear operator on $H(X)$ defined by $M_\phi f=\phi f$.  A set $\Omega\subset X$ is a \textbf{set of uniqueness} 
for $Mult(X)$ if $\phi_1=\phi_2$ in $Mult(H(X))$ whenever $\phi_1,\phi_2\in Mult(H(X))$ and $\phi_1=\phi_2$ on $\Omega$. 
We use $Mult_1(H(X))$ to denote the unit ball of $Mult(H(X))$. 

The Nevanlinna-Pick problem on $H(X)$ is to determine, given nodes $\l_1,...,\l_N\in X$ and 
target points $\oo_1,...,\oo_N\in\D$, whether there exists a $\phi\in Mult_1(H(X))$ that satisfies $\phi(\li)=\oo_i$ for each $i$. For each $\l\in X$, we use 
$k_{\l}$ to denote the reproducing kernel at $\l$, the element of $H(X)$ that represents the linear functional evaluation at $\l$. Given a problem on 
$H(X)$ with nodes $\l_1,...,\l_N\in X$ and target points $\oo_1,...,\oo_N\in\D$, let $k_{ij}=<k_{\l_j},k_{\li}>$. We define two $N$ by $N$ matrices $W$ and 
$K$ with the following formulas,
\[W = ( 1 - \oo_i \bar \oo_j) \, \textrm{ and } \, K=(k_{ij} ) .\]
Finally, we use $W\cdot K=((1-\oo_i \bar\oo_j)k_{ij})$ to denote the Schur entrywise product of $W$ and $K$, and call $W\cdot K$ the Pick matrix associated to the problem.

\bd
\label{Pick-like}
Consider a Hilbert function space $H(X)$ with kernel $K$. We say that $K$ is a \textbf{Pick uniqueness kernel} if the following holds for each solvable 
Nevanlinna-Pick problem on $H(X)$: If the Pick matrix $W\cdot K$ associated to the problem is singular, then the problem has a unique solution.
\ed

\bt
\label{FunctionSpace}
Consider a Hilbert function space $H(X)$ with kernel $K$. If for each finite set of points $\l_1,...,\l_N\in X$ and non-zero scalars $a_1,...,a_N$ the complement of the zero set of $G(x)=a_1 k_{\l_1}(x)+...+a_N k_{\l_N}(x)$ is a set of uniqueness for $Mult(H(X))$, then $K$ is a Pick uniqueness kernel.
\et

We prove Theorem \ref{FunctionSpace} in section 2, by generalizing an argument from \cite{agmc_dv}.


\subsection{Theorem \ref{ThmB}}
Given 1-dimensional inner variety $V\subset\C^n$, Theorem \ref{ThmB} gives a geometric procedure for constructing a problem with the property that 
$V\cap\D^n\subset U$. We will use Theorem \ref{ThmB} us to establish one half of Theorem \ref{EachInner}.

A rational function $F=\frac{r}{q}$ with $q,r$ relatively prime is called \textbf{inner} if $q\neq0$ on $\D^n$ and $|F(\tau)|=1$ for almost every $\tau\in\T^n$. A rational inner function $F$ is called \textbf{regular} if $q\neq 0$ on $\overline{\D^n}$. Given a 1-dimensional inner variety $V\subset\C^n$, the rank of $V$ is the n-tuple of the generic number of sheets of $V$ over each coordinate. We define the \textbf{degree of $F$ on $V$}, $\deg_V(F)$, as the valence of $F$ on 
$V\cap\D^n$.

\bt
\label{ThmB}
Consider a regular rational inner function $F$ on $\D^n$ and an irreducible 1-dimensional inner variety  $V$. If $\deg_V(F)<N$ and $\l_1,...,\l_N\in V\cap \D^n$ 
are distinct, then every solution to the Nevanlinna-Pick problem with data $\l_1,...,\l_N$ and $F(\l_1),...,F(\l_N)$ equals $F$ on $V\cap\D^n$, i.e. $V\cap\D^n\subset U$.
\et



We prove Theorem \ref{ThmB} in section 3. We give the formula for $\deg_V(F)$, from \cite{agmc_dv}, in Theorem \ref{FredholmDegree} of section 3 .

\subsection{Theorem \ref{ThmD}}
Given an inner variety $V\subset\C^n$ and a problem with nodes on $V$, Theorem \ref{ThmD} complements Theorem \ref{ThmB} by establishing sufficient conditions for $U\subset V\cap\D^n$. We will use Theorem \ref{ThmD} to prove Theorems \ref{EachInner} and \ref{kUniqueSet}, and to partially answer a question by Agler and McCarthy from \cite{agmc_dv}.

Consider a rational inner function $F=\frac{r}{q}$ on $\D^n$ with $q,r$ relatively prime. We define the \textbf{n-degree} of $F$, $n$-$\deg(F)$, as the $n$-tuple with $i$th entry given by the degree of $r$ in $z_i$. Given a polynomial $p$ and a rational inner function $F$, we use $Z_p$ to denote the zero set of $p$ and write $n$-$\deg(p)\leq n$-$\deg(F)$ if the degree of $p$ in each variable $z_i$ is less than or equal to the degree of $r$ in $z_i$. For an inner variety $V\subset\C^n$, we use the notation $n$-$\deg(V)\leq n$-$\deg(F)$ to mean that there exist irreducible polynomials $p_1,...,p_r$ such that 
$V=Z_{p_1}\cap...\cap Z_{p_r}$ and $n$-$\deg(p_k)\leq n$-$\deg(F)$ for each $k$. If a rational inner function $F_\epsilon$ satisfies $F_\epsilon=F$ on 
$V\cap\D^n$ and $F_\epsilon\neq F$ on $\D^n$, then we call $F_\epsilon$ a \textbf{rational inner perturbation of $F$ on $V$}. 

\bt
\label{ThmD}
If $V$ is  an inner variety and $F$ is a regular rational inner function on $\D^n$ with $n$-$\deg(V)\leq n$-$\deg(F)$, then
for every $z \in \D^n\setminus V$ there exists a rational inner perturbation of $F$ on $V$, $F_\epsilon$, such that $F_\epsilon(z)\neq F(z)$.
\et

We prove Theorem \ref{ThmD} in Section 4. The following corollary of Theorem \ref{ThmD} provides a partial answer to Agler and McCarthy's Question 4.14 in \cite{agmc_dv}. Further results related to their question will appear in \cite{dsQ}.


\begin{corollary}
\label{PartialAnswer}
Consider an extremal Nevanlinna-Pick problem on $\D^2$ and an inner variety $V$ given as the zero set of a polynomial $p$. If there exists a regular rational inner function $F$ solving the problem such that the degree of $F$ is greater than or equal to the degree of $p$ in each variable, then $V\cap\D^2=U$.
\end{corollary}

\subsection{Organization of the Paper}
The rest of this paper is organized as follows. In section 2, we prove Theorem \ref{FunctionSpace}. 
In section 3, we show that for each 1-dimensional inner variety $V$, there exists a Hilbert function space on $V\cap\D^n$ the kernel for which is a Pick uniqueness kernel. In section 4 we prove Theorem \ref{ThmB}. In section 5 we prove Theorem \ref{ThmD}. In section 6, we prove Theorems \ref{EachInner} and \ref{kUniqueSet}. In section 7 we use Theorems \ref{ThmB} and \ref{ThmD} to construct a Nevanlinna-Pick problem on $\D^2$ for which the set of uniqueness equals the Neil Parabola. 

I would like to thank Jim Agler for his generous help in discussing the ideas leading to this work and for his help in vastly improving the exposition of this paper. 
In particular, the idea for Lemma \ref{InnerRank} was privately communicated to the author by Jim Agler and inspired much of this work.


\section{Proof of Theorem 1.2}
This section is dedicated to the proof of Theorem \ref{FunctionSpace}. Accordingly, let $H(X)$ be a Hilbert function space with kernel $K$ and the property that for 
each finite set of points $\l_1,...,\l_N\in X$ and non-zero scalars $a_1,...,a_N$ the complement of the zero set of $G(x)=a_1 k_{\l_1}(x)+...+a_N k_{\l_N}(x)$ is a set of uniqueness for $Mult(H(X))$. Fix a Nevanlinna-Pick problem with data $\l_1,...,\l_N\in X$ and $\oo_1,...,\oo_N\in\D$, let $k_{ij}=<k_{\l_j},k_{\l_i}>$. Let 
$W$ and $K$ denote the following $N$ by $N$ matrices 
\[W = ( 1 - \oo_i \bar \oo_j) \, \textrm{ and }  \, K=(k_{ij} ).\]
Assume that $W\cdot K=((1-\oo_i \bar\oo_j)k_{ij})$ is singular and that $\phi_1,\phi_2\in Mult_1(H(X))$ are solutions to the problem. We seek to show that $\phi_1=\phi_2$.
Let $\gamma$ be a non-zero vector in the null space of $W\cdot K$, let $\l_{N+1}$ be any point in $X$ that is not one of the original nodes and 
let $\wN$ be a possible value that a solution to the original problem can take at $\lN$. Since the problem with data $\l_1,..,\l_{N+1}$ and 
$\oo_1,...,\oo_{N+1}$ is solvable, Theorem 5.2 of \cite{ampi} implies that the matrix $[ (1 - w_i \overline{w_j}) k_{ij}]_1^{N+1}$ is positive semi-definite 
(one can check this by computing the operator $1-M_\phi M_\phi^*$ on the span of $\{k_{\li}\}$). 
Thus, for each $z \in \C$,
\be
\label{eqd6}
\la \left[ (1 - w_i \overline{w_j}) k_{ij} \right]_1^{N+1} 
\left( \begin{array}{c}
\gamma\\
z \end{array} \right) \ , \ 
\left( \begin{array}{c}
\gamma\\
z \end{array} \right)\ra \quad \geq \quad 0.
\ee
Since $\gamma$ is in the null-space of $[ (1 - w_i \overline{w_j})k_{ij}]_1^N$, inequality~(\ref{eqd6}) reduces to 
\be
\label{eqd7}
2 \Re [ \bar z \sum_{j=1}^N (1 - \bar w_j \wN) \g_j k_{N+1,j} ]
\ + \ |z|^2 (1 - |\wN|^2) ||k_{\l_{N+1}}||^2 \ \geq \ 0 .
\ee
Since equation \ref{eqd7} holds for all $z$, it follows that  
\[\sum_{j=1}^N (1 - \bar w_j \wN) \g_j k_{N+1,j}=0\]
and the following gives an implicit formula for $\wN$,

\be
\label{eqd799}
\wN \ \left( \sum_{j=1}^N \bar w_j \g_j k_{N+1,j} \right) = \sum_{j=1}^N  \g_j k_{N+1,j}.
\ee

\noindent
Define $G\in H(X)$ with the formula 
\[G(\l_{N+1}) = \sum_{j=1}^N \g_j k_{\l_j} (\l_{N+1})= \sum_{j=1}^N \g_j k_{N+1,j}.  \] 
Let $\Omega\subset X$ be the complement of the zero set of $G$, i.e. the set of points $\l_{N+1}\in X$ on which both sides of formula \ref{eqd799} do not vanish. 
On $\Omega$, the value of $\wN$ is uniquely determined by the following formula,
\be
\label{eqd8}
\wN = \left(\sum_{j=1}^N  k_{N+1,j} \g_j \right)/ \left(\sum_{j=1}^N \bar w_j k_{N+1,j} \g_j\right).
\ee
\noindent
Thus, $\phi_1=\phi_2$ on $\Omega$ and the assumption that $\Omega$ is a set of uniqueness for $Mult(H(X))$ implies that $\phi_1=\phi_2$. \ep


\section{A Hilbert function space on $V$}
In this section we prove the following theorem and two related results. The Schur Class of $\D^n$, $\Scn$, is the set of analytic functions that map $\D^n$ to 
$\overline{\D}$, i.e. satisfy $||F||_\infty\leq1$.  

\bt
\label{measure}
Given a 1-dimensional irreducible inner variety $V\subset\C^n$, there exists a finite measure $\mu$ on $V\cap\T^n$ such that $\htm$, the closure of the polynomials in $L^2(\mu)$, has the following properties.\\
\textbf{i.} The set of kernel functions $\{k_\l:\l\in V\cap\D^n\}$ is dense in $\htm$.\\
\textbf{ii.} The kernel associated to $\htm$ is a Pick uniqueness kernel.\\
\textbf{iii.} For each $F\in \Scn$, if $f$ is the restriction of $F$ to $V\cap\D^n$, then\\
$f\in Mult(\htm)$ and $||M_f||\leq||F||_\infty$.
\et

Theorem \ref{measure} was inspired by a similar result for $n=2$ in \cite{agmc_dv}.

\noindent
\bp
Let $p_1,...,p_r$ be irreducible polynomials such that $V=Z_{p_1}\cap...\cap Z_{p_r}$. Let $C$ be the projective closure of 
$Z_{p_1}\cap...\cap Z_{p_r}$ in $\C {\mathbb P}^n$ and identify $V\cap\D^n$ with a subset of $C$. Let $(S,\phi)$ be the desingularization of $C$, 
a compact Riemann surface $S$ and a holomorphic function $\phi : S \to C$ that is biholomorphic from $S'$ onto $C'$ and finite-to-one from $S \setminus S'$ onto 
$C \setminus C'$. Here $C'$ is the set of non-singular points in $C$, and $S'$ is the preimage of $C'$. See e.g. \cite{gh78} for details of the desingularization.


Let $\Om = \phi^{-1} ( V\cap\D^n)$. Then $\partial \Om$ is a finite union of disjoint curves, each of which is analytic except possibly at a finite number of cusps. 
Lemma 1.1 of \cite{agmc_dv} states that there exists a finite measure $\nu$ on $\partial \Om$ such that evaluation at every $\l$ in $\Om$ is a bounded linear functional on $A^2(\nu)$, the closure in $L^2(\nu)$ of $A(\Om)$, the functions that are analytic on $\Omega$ and continuous on $\overline{\Omega}$. Furthermore, 
the lemma states that the linear span of the corresponding evaluation kernels is dense in $A^2(\nu)$.


The desired measure $\mu$ is the push-forward of $\nu$ by $\phi$, normalized to have mass 1 on $\partial V:=\partial (V\cap\D^n)=V\cap\T^n$.
In particular, $\mu$ is defined by
\[\mu(E) = \nu (\phi^{-1}(E)).\]
If $\l\in V\cap\D^n$ is a regular point of $V$, then there exists a unique $\z$ such that $\phi(\z) = \l$ and $k_\z$ is the reproducing kernel associated to 
$\z$ in $A^2(\Om)$. 
The function $k_\l=k_\z \circ \phi^{-1}$ is defined $\mu$ almost everywhere on $\partial V$ and for each $f\in\htm$ satisfies
\[
<f,k_\l>_{\htm}=\int_{\partial V} f \cdot \overline{k_\z \circ \phi^{-1}}  d \mu
\=\
\int_{\partial \Om} f \circ \phi  \cdot \overline{k_\z} d\nu
\= <f\circ\phi,k_\z>_{A^2(\Om)} \= f(\l).
\]
If $\l\in V\cap\D^n$ is a singular point of $V$, then there exist finitely many $\z_1,...,\z_s$ such that $\phi(\z_i) = \l$ and the function
\[k_\l=\frac{1}{s}\sum k_{\z_i} \circ \phi^{-1}\]
is the corresponding reproducing kernel function for $\l$.

To see that the kernel associated to $\htm$ is a Pick uniqueness kernel fix points $\l_1,...,\l_N\in V\cap\D^n$, fix non-zero scalars $\overline{a_1},...,\overline{a_N}$, 
fix the function $G(x)=a_1 k_{\l_1}(x)+...+a_N k_{\l_N}(x)$ and let $O\subset V\cap\D^n$ be the complement of the zero set of $G(x)$. 
That $O$ is a set of uniqueness for $Mult(\htm)$ follows immediately from the following two observations.\\

\noindent
\textbf{1:} $Mult(\htm)\subset\htm$. This holds since $1\in\htm$.\\
\textbf{2:} Each function $f\in\htm$ that vanishes on an open $O\subset V\cap\D^n$ vanishes on $V\cap\D^n$. This holds since the polynomials are dense in 
$\htm$ and if a polynomial $p$ vanishes on $O\subset V\cap\D^n$, then the function $P=p\circ \phi$ is analytic on $\Om$ and vanishes on $\phi^{-1}(O)$ and thus, $p$ vanishes on $V\cap\D^n$.


To establish part \textbf{iii.} of the theorem, fix $F\in \Scn$ and let $f$ be the restriction of $F$ to $V\cap\D^n$. Since $||F||_\infty\leq1$, there exist polynomials 
$\{p_i\}$ that satisfy $||p_i||_\infty\leq1$ and approximate $F$ on $\D^n$ and in particular on $V\cap\D^n$. Thus, $f\in \htm$. The following calculation shows 
that $||M_f||\leq 1$. For each $g\in\htm$ we have that\\
$\displaystyle{ ||M_f g||=\int_{\partial V} |M_f g|^2  d \mu=\int_{\partial V} |f g|^2  d \mu\leq\int_{\partial V} |g|^2  d \mu=||g|| }.$ \ep

In Theorem \ref{FredholmDegree} and Lemma \ref{InnerRank} below, we establish two useful properties of $\htm$. To state them, we need some notation.
 
Given a polynomial $q(z_1,...,z_n)$ with $n$-$\deg(q)=(d_1,...,d_n)$, let
\be z^d= z_1^{d_1}\cdot\cdot\cdot z_n^{d_n}, \hspace{10pt} \frac{1}{\overline{z}}=( \frac{1}{\overline{z_1}},...,\frac{1}{\overline{z_n} }) \hspace{10pt} \textrm{ and } \hspace{10pt} 
\widetilde{q}(z)=z^d\overline{q(\frac{1}{\overline{z}})}. 
\ee

\bt(Rudin \cite{rud69})
\label{Rudin}
Given a rational inner function $F$ on $\D^n$, there exist a polynomial $q$ that does not vanish on $\D^n$, an n-tuple of positive integers $m=(m_1,...,m_n)$ 
and $\tau\in\T$  such that
\be
F(z)=\tau z^m\frac{\widetilde{q}}{q}.
\label{InnerFormula}
\ee
Furthermore, each rational function $F$ of the form \ref{InnerFormula} is inner. 
\et

For a Hilbert space $H$, a bounded linear operator $A$ on $H$ is a \textbf{Fredholm operator} if it has closed range, $dim(Ker(A)) < \infty$ and $dim(ker(A^*)) < \infty$. If $A$ is a Fredholm operator, then the Fredholm index of $A$ is defined to be $ind(A) = dim(Ker(A))-dim(Ker(A))$.
The following theorem summarizes the well known results that we will use.
\bt{(Conway, Ch 12. \cite{Conway90})} .\\
\label{Fredholm}
\textbf{i.} For $A,B$ Fredholm operators on $H$, $ind(AB)=ind(A)+ind(B)$.\\
\textbf{ii.} The set of Fredholm operators, denoted $\mathcal{F}$, is open in $B(H)$.\\
\textbf{iii.} As a function, $ind:\mathcal{F} \to \mathbb{Z}$, is constant on connected components of $\mathcal{F}$.
\et

Given a 1-dimensional inner variety $V\subset \C^n$ and a rational inner function $F$ on $\D^n$, for almost every $z\in\D$, the cardinality of the set 
$\{\l\in V\cap\D^n:F(\l)=z\}$ constant. We define the \textbf{degree of $F$ on $V$} by fixing one such $z\in\D$ and letting 
$\deg_V(F)=|\{\l\in V\cap\D^n:F(\l)=z\}|$. We modify a proof from \cite{agmc_dv} to establish the following result.

\bt
\label{FredholmDegree}
Consider a 1-dimensional inner variety $V\subset \C^n$ with $rank(V)=(m_1,...,m_n)$, $\htm$ the Hilbert function space from Theorem \ref{measure}, 
$F$ a regular rational inner function with $n$-$deg(F)=(d_1,...,d_n)$ and $M_F^*$ the adjoint of the bounded linear operator of multiplication by $F$ on $\htm$. Then,
\[ \deg_V(F) =m_1d_1+...+m_n d_n=ind(M_F^*) \]
\et

\bp
We first establish the following for each coordinate function $z_i$, 
\[\deg_V(z_i)=m_i =m_i-0=dim(ker(M_{z_i}^*))-dim(ker(M_{z_i}))=ind(M_{z_i}^*).\] 
After an automorphism of $\D^n$, we may assume that $V$ has $m_i$ distinct non-singular points $\l_{i,1},...,\l_{i,m_i}$ lying over the zero set of each coordinate 
$z_i$. That $\deg_V(z_i)=m_i$ is immediate. The operator of multiplication by $z_i$, denoted $M_{z_i}$, is an isometry with finite multiplicity and thus, 
Fredholm. 

Let $K=kernel(M_{z_i}^*)$ and let $K_i=span\{k_{\l_{i,1}},...,k_{\l_{i,m_i}}\}$. Since kernel functions are eigenvectors of the adjoints of multiplication operators 
and satisfy $M_F^*k_{\l}=\overline{F(\l)}k_{\l}$, we have that $K_i\subset K$. Suppose towards a contradiction that $K_i\subsetneq K$. Choose $g\in K\ominus K_i$. 
That $g$ is orthogonal to $K_i$ means that $g$ vanishes at each of the points $\l_{i,1},...,\l_{i,m_i}$ and that the function $\frac{g}{z_i}$ is in $\htm$. Thus, 
$g$ equals $M_{z_i}\frac{f}{z_i}$ and thus, $g\in range(M_{z_i})\subset H\ominus K$, a contradiction.

For a regular rational inner function $F(z)=\frac{\widetilde{q}}{q}$, normalize $q$ so that $q(0)=1$, let $q(z)=1+Q(z)$ and for $0\leq r\leq1$ let
\[ F_r(z)=\frac{z^d+z^d\overline{Q(r\frac{1}{\overline{z}})}}{1+Q(rz)}\]
For $F_0=z^d$, the first part of the following equality is immediate and the second part follows from part \textbf{i} of Theorem \ref{Fredholm}.
\[\deg_V(F_0)=m_1d_1+...+m_nd_n=ind(M_{F_0}^*)\]
For each $r$ between $0$ and $1$, $F_r$ is a regular rational inner function and thus, $F_r(V\cap\T^n)\subset \T$ implying that $\deg_V(F_r)$ remains constant and 
$M_{F_r}^*$ is a Fredholm operator meaning that part \textbf{iii} of Theorem \ref{Fredholm} implies that $ind(M_{F_r}^*)$ remains constant.
\ep

We can now prove another useful property of $\htm$.

\bl
\label{InnerRank}
Consider an irreducible 1-dimensional inner variety $V\subset\C^n$, a rational inner function $F$ on $\D^n$ and  
$\htm$ the Hilbert function space on $V\cap\D^n$ from Lemma \ref{measure}. Let $\l_1,...,\l_N\in V$ be distinct, let $\oo_i=F(\li)$ and let $k_{ij}=<k_{\l_j},k_{\l_i}>$. 
The rank of the $N$ by $N$ matrix $W\cdot K=\displaystyle{((1 - w_i \overline{w_j}) k_{ij})}$ is less than or equal to $\deg_V(F)$.
\el

\bp
Fix a measure $\mu$ on $V\cap\T^n$, the existence of which is guaranteed by lemma \ref{measure}, and consider multiplication by $F$, denoted 
$M_F$, as a bounded linear operator on $\htm$. The Pick matrix $W\cdot K$ equals the grammian of the vectors $\{(1-M_F M_F^*)k_{\li}\}_1^N$, i.e.
\[ <(I- M_F  M_F^*)k_{\lj},(I- M_F  M_F^*)k_{\li}> =<(I-M_F  M_F^*)k_{\lj},k_{\li}> = \]
\[<k_{\lj},k_{\li}>-< M_F^* k_{\lj}, M_F^* k_{\li}>=<k_{\lj},k_{\li}>-<\overline{F(\lj)}k_{\lj},\overline{F(\li)}k_{\li}> = \]
\[<k_{\lj},k_{\li}>-\oo_i\overline{\oo_j}<k_{\lj},k_{\li}>=(1-\oo_i\overline{\oo_j})K_{i,j}\]

The conclusion of the theorem follows from the following inequality,
\be \label{rank} rank(W\cdot K)\leq\dim (sp\{(1-M_F M_F^*)k_{\li}\}) \leq \deg_V(F). \ee

The first inequality in \ref{rank} follows from the fact that the rank of the grammian of a set of vectors is less than or equal to the dimension of the span of the vectors. 
The second inequality in \ref{rank} follows from the observation that $range(1-M_F M_F^*)=kernel(M_F^*)$ and the equality $dim(kernel(M_F^*))=\deg_V(F)$ from Theorem \ref{FredholmDegree}.
\ep

\section{Proof of Theorem \ref{ThmB}}
Lemmas \ref{measure} and \ref{InnerRank} of the previous section readily imply Theorem \ref{ThmB}.

\noindent
{\em Proof of Theorem \ref{ThmB}: } Fix a 1-dimensional inner variety $V\subset \C^n$, fix a regular rational inner function $F$ on $\D^n$, fix $N>\deg_V(F)$, fix 
distinct points $\l_1,..,\l_N\in V\cap\D^n$ and let $\oo_i=F(\li)$. Suppose that $G\in\Scn$ satisfies $G(\li)=F(\li)$ for each $i=1,...,N$. 
We seek to show that $G=F$ on $V\cap\D^n$.

Let $\htm$ be the Hilbert function space from Lemma \ref{measure}. Theorem \ref{measure} states that the kernel associated with $H$ is a Pick uniqueness kernel. 
Lemma \ref{InnerRank} states the $N$ by $N$ Pick matrix $W\cdot K$ corresponding to Nevanlinna-Pick problem with data $\l_1,..,\l_N$ and $\oo_1,...,\oo_N$ has 
rank less than or equal to $\deg_V(F)$ and is thus, singular. It follows that $F$ is the unique solution to the problem in $Mult_1(H)$ and since $G$ is another solution, we have that $G=F$ in $Mult(\htm)$ and on $V\cap\D^n$. \ep


\section{Proof of Theorem \ref{ThmD}}




This section is dedicated to the proof of Theorem \ref{ThmD}. Accordingly, let $V$ be an inner variety and fix $F$ a regular rational inner function on $\D^n$ with 
$n$-$\deg(V)\leq n$-$\deg(F)$. Let $F=z^{m}\frac{\widetilde{q}}{q}$, let $s=n$-$\deg(q)$, let $p_1,...,p_k$ be irreducible polynomials such that 
$V=Z_{p_1}\cap...\cap Z_{p_k}$ and let $r_i=n$-$\deg(p_i)$. 

For each pair of $k$-tuples of real numbers $\epsilon$ and $\delta$ define the function $F_{\epsilon,\delta}$ with the following formula,
\[
\label{Perturb3}
F_{\epsilon,\delta}=\Frac{z^{m+s}\left(\overline{q(\frac{1}{\bar z})} +\sum_1^r \epsilon_i \overline{p(\frac{1}{\bar z})} _i + 
\delta_i \overline{\widetilde{p}_i(\frac{1}{\bar z}) }\right) }{q+\sum_1^r \epsilon_i p_i + \delta_i \widetilde{p_i}}.
\]
The assumption that $F$ is regular means that $q$ does not vanish on $\overline{\D^n}$ and thus, when the entries of $\epsilon$ and $\delta$ are sufficiently small, Theorem \ref{Rudin} implies that the function $F_{\epsilon,\delta}$ is inner. The assumption that $n$-$\deg(V)\leq n$-$\deg(F)$ 
implies that for each $i$ we have that each entry of the $n$-tuple $m+s-r_i$ is non-negative and thus,
\[
\label{Perturb3}
F_{\epsilon,\delta}=\Frac{z^m\widetilde{q}+\sum_1^r \epsilon_i z^{m+s-r_i}\widetilde{p_i} + \delta_i z^{m+s-r_i}p_i}{q+\sum_1^r \epsilon_i p_i + \delta_i \widetilde{p_i}}.
\]
Since $V$ is inner, Proposition 2.6 of \cite{ams06} implies for each $i$ the polynomial $\widetilde{p_i}$ vanishes on $V$ and thus, 
$F_{\epsilon,\delta}=z^{m}\frac{\widetilde{q}}{q}=F$ on $V$.

To conclude the proof of the theorem, fix a $z\in\D^n \smallsetminus V$ and choose a $p_h$ such that $p_h(z)\neq0$. Fixing $\epsilon_h$ and letting each 
$\delta_i$ and each $\epsilon_i$ with $i\neq h$ go to zero results in a function that satisfies $F_{\epsilon,\delta}(z)\neq F(z)$.
\ep

The case where $F=\frac{\widetilde{q}}{q}$, the variety $V=Z_p$ is the zero set of a single polynomial $p$, and $n$-$\deg(p)=n$-$\deg(F)$ is particularly simple. 
For $\epsilon$ sufficiently small, the function given by the  following formula is a rational inner perturbation of $F$ on $V$.
\be
\label{Perturb0}
F_\epsilon=\Frac{\widetilde{q+\epsilon p}}{q+\epsilon p}=\Frac{\widetilde{q}+\epsilon\widetilde{p}}{q+\epsilon p}=\Frac{\widetilde{q}+\epsilon p}{q+\epsilon p}.
\ee
\section{Proof of Theorems \ref{EachInner} and \ref{kUniqueSet}}
Theorems \ref{ThmB} and \ref{ThmD} readily imply Theorem \ref{EachInner}.

{\sc proof of Theorem \ref{EachInner}:} Fix a 1-dimensional inner variety $V$, let $F$ be a regular rational inner function $F$ on $\D^n$ such that 
$n-\deg(V)\leq n-\deg(F)$, let $N>\deg_V(F)$ and fix distinct points $\l_1,...,\l_N\in V\cap\D^n$. By Theorem \ref{ThmB}, the set of uniqueness 
$U$ for the Nevanlinna-Pick problem with data $\l_1,...,\l_N$ and $F(\l_1),...,F(\l_N)$ satisfies  $V\cap\D^n \subset U$. For each $z \not\in V\cap\D^n$, Theorem 
\ref{ThmD} guarantees the existence of an $F_\epsilon$ that solves the problem and satisfies $F_\epsilon(z)\neq F(z)$. Thus, $z\not \in U$ and $U\subset V\cap\D^n$.
\ep

Theorem \ref{ThmD} and Theorem 1.1 from \cite{dsPLn} readily imply Theorem \ref{kUniqueSet}.

{\sc Proof of Theorem~\ref{kUniqueSet}:}
Fix $n\geq2$, $k\leq n$. Fix $n$ regular rational inner functions $m_1,...,m_n$ on $\D^k$, non of which are constant in any the variables $z_i$, and define the analytic mapping $D:\D^k\to\D^n$ with the following formula 
\[D(z)=(m_1(z),...,m_n(z)).\] 
It follows that $D(\T^k)\subset\T^n$ and that there exists a $k$-dimensional inner variety $V$ such that $V\cap\D^n=D(\D^k)$. 
Fix an analytic function $F$ on $\D^n$ such that $n$-$\deg(V)\leq n$-$\deg(F)$. Identify the restriction of $F$ to $V\cap\D^n$ with a rational inner function $f$ on 
$\D^k$, $f(z_1,...,z_k)=F(D(z_1,...,z_k))$. Theorem 1.1 of \cite{dsPLn} guarantees the existence of a Nevanlinna-Pick problem with data $\l_1,...,\l_M\in \D^k$ and $f(\l_1),...,f(\l_M)$ for which the solution is unique on $\D^k$. Consider the Nevanlinna-Pick problem on $\D^n$ with data $D(\l_1),...,D(\l_M)\in D^n$ and 
$f(\l_1),...,f(\l_M)$. If $G$ is a solution to the new problem on $\D^n$, then $g(z_1,...,z_k)=G(D(z_1,...,z_k))$ is a solution to the original problem on $\D^k$ and thus, must equal $f$ on $\D^k$. Thus, $G=F$ on $V\cap\D^n$ and it follows that $V$ is a subset of the set of uniqueness for the problem on $\D^n$, i.e. $V\cap\D^n\subset U$. 
Since $n$-$\deg(V)\leq n$-$\deg(F)$, Theorem \ref{ThmD} implies that $U\subset V\cap\D^n$.
\ep

\section{The Neil Parabola}

Let $\mathcal{N}$ denote the Neil Parabola, the zero set of the polynomial $p=z^3-w^2$. Below, we use Theorems \ref{ThmB} and \ref{ThmD} can be used to construct a Nevanlinna-Pick problem on $\D^n$ for which the set of uniqueness equals $\mathcal{N}\cap\D^2$.

Let $F(z,w)=z^3w^2$, let $N>12=\deg_\mathcal{N}(F)$ and fix a Nevanlinna-Pick problem with data $\l_1,...,\l_N\in V$ and $F(\l_1),...,F(\l_N)$. Theorem \ref{ThmB} states that all solutions to this problem agree on $V$, i.e. $\mathcal{N}\cap\D^2\subset U$. Since  $n-\deg(\mathcal{N}) = (3,2)\leq (3,2) =n-\deg(F)$, Theorem \ref{ThmD} implies that $U\subset\mathcal{N}\cap\D^2$. 
That $U\subset\mathcal{N}\cap\D^2$ also follows directly from the fact that for each $\epsilon<1/2$ the function $F_\epsilon$, given by the formula 
below, is a rational inner perturbation of $F$ on $V$. 
\[ F_{\epsilon}(z,w)=\frac{z^3w^2+\epsilon \widetilde{p}}{1+\epsilon p}=\frac{z^3w^2+\epsilon(w^2-z^3)}{1+\epsilon (z^3-w^2)}=\frac{z^3w^2-\epsilon p}{1+\epsilon p}. \]

\end{document}